\documentclass[reqno]{amsart}
\usepackage{hyperref}

\newcommand{\be}{\begin{equation}}
\newcommand{\ee}{\end{equation}}
\newcommand{\bea}{\begin{eqnarray}}
\newcommand{\eea}{\end{eqnarray}}
\newcommand{\bean}{\begin{eqnarray*}}
\newcommand{\eean}{\end{eqnarray*}}
\newcommand{\la}{\label}
\newcommand{\R}{{\mathbb{R}}}
\newcommand{\C}{{\mathbb{C}}} 
\newcommand{\re}{{\rm Re}\;}
\newcommand{\cE}{{\mathcal E}}
\newcommand{\cL}{{\mathcal L}}
\newcommand{\lf}{\lambda\phi}
\newcommand{\nl}{\newline}
\newcommand{\po}{\partial\Omega}
\newcommand{\tarr}[6]{{\left\{\begin{array}{lll}
{#1}&{#2}\\[0.2cm]
{#3}&{#4}\\[0.2cm]
{#5}&{#6}
\end{array}\right.}}

\begin{document}
\title[Heat kernel estimates]
{Heat kernel estimates for fourth order non-uniformly elliptic operators with non strongly convex symbols}

\author{Gerassimos Barbatis, Panagiotis Branikas}

\address{Gerassimos Barbatis \newline
Department of Mathematics,
National and Kapodistrian University of Athens, \nl
Panepistimioupolis, 15784 Athens, Greece}
\email{gbarbatis@math.uoa.gr}

\address{Panagiotis Branikas \newline
Department of Mathematics,
National and Kapodistrian University of Athens, \nl
Panepistimioupolis, 15784 Athens, Greece}
\email{pbranikas@math.uoa.gr}


\begin{abstract}
We obtain heat kernel estimates for a class of fourth order non-uniformly elliptic operators in two dimensions. 
Contrary to existing results, the operators considered have symbols that are not strongly convex. This entails certain
difficulties as it is known that, as opposed to the strongly convex case, there is no absolute exponential constant. Our estimates involve sharp constants and Finsler-type distances that are induced by the operator
symbol. The main result is based on two general hypotheses, a weighted Sobolev inequalitry and an interpolation inequality, which are related to the singularity or degeneracy of the coefficients.
\end{abstract}

\maketitle
\numberwithin{equation}{section}
\newtheorem{theorem}{Theorem}[section]
\newtheorem{lemma}[theorem]{Lemma}
\newtheorem{proposition}[theorem]{Proposition}
\newtheorem{remark}[theorem]{Remark}
\allowdisplaybreaks

\section{Introduction}

Let $\Omega$ be a planar domain and let
\be
Hu = \partial_{x_1}^2 \big( \alpha(x) \partial_{x_1}^2 u  \big) +2\partial_{x_1x_2}^2 \big( \beta(x) \partial_{x_1x_2}^2 u  \big)
+ \partial_{x_2}^2 \big( \gamma(x) \partial_{x_2}^2 u  \big)
\label{oper h}
\ee
be a fourth-order,  self-adjoint, uniformly elliptic operator in divergence form on $\Omega$ with measurable coefficients
satisfying Dirichlet boundary conditions on $\partial\Omega$. It has been shown by Davies \cite{d} that $H$ has a continuous heat kernel
$G(x,x',t)$ which satisfies the Gaussian-type estimate
\be
|G(x,x',t)|\leq c_1 t^{-\frac{1}{2}}\exp\Big(-c_2\frac{|x-x'|^{4/3}}{t^{1/3}}+c_3t\Big),
\la{eq:2}
\ee
for some positive constants $c_1$, $c_2$, $c_3$ and all $t>0$ and $x,x'\in\Omega$. Indeed \cite{d} deals with the more general case of an operator
of order $2m$ acting on a domain in $\R^n$,
$n<2m$.
The study of fundamental solutions is central in the theory
of linear parabolic PDEs. For more results on
heat kernel estimates for higher-order operators we refer to
\cite{CLYZ,Cl,DDY,Du,DzHe,HWZD,RS-C1,RS-C2}.
See also \cite{QuR-B,Ze} for related results specific to fourth order
operators.

A sharp version of the Gaussian estimate (\ref{eq:2}) is obtained  in \cite{b2001} where it was proved that 
\be
|G(x,x',t)|\leq c_{\epsilon} t^{-\frac{1}{2}}\exp\Big\{-\Big( \frac{3\sqrt[3]{2}}{16} -c\theta  -\epsilon\Big) \frac{d_M(x,x')^{4/3}}{t^{1/3}}+
c_{\epsilon,M}t\Big\},
\la{eq:3}
\ee
for arbitrary $\epsilon$ and $M$ positive. Here $\theta \geq 0$ is a constant that is related to the regularity of the coefficients and $d_M(x,x')$, $M>0$, is a family of Finsler-type distances on $\Omega$ which
is monotone increasing and converges as $M\to +\infty$ to a limit Finsler distance $d(x,x')$. The sharpness follows by comparing against the short time asymptotics obtained in
\cite{ep} for equations with constant coefficients and which involve precisely the constant $3\sqrt[3]{2}/16$;
we refer to \cite{bb2018} for a more detailed discussion
of the distance function $d(x,x')$.

An important assumption for both the Gaussian esimate (\ref{eq:3}) and for the
corresponding asymptotic estimate of \cite{ep}
is the {\em strong convexity} of the symbol
\be
A(x,\xi) =\alpha(x) \xi_1^4 +2\beta(x) \xi_1^2 \xi_2^2 +\gamma(x)\xi_2^4 \; , \qquad x\in\Omega \,  ,  \;\; \xi \in\R^2 \, ,
\la{symbol}
\ee
of the operator $H$. The notion of strong convexity was introduced in \cite{ep} and it applies to operators 
of order $2m$ acting on $\R^d$ which have constant coefficients.
In our context the strong convexity of the symbol (\ref{symbol})  amounts to
\be
0\leq  \beta(x) \leq 3 \sqrt{\alpha(x)\gamma(x)} \, , \qquad x\in \Omega \, .
\footnote{We note that (\ref{str conv}) is the assumption made for the heat kernel estimates of \cite{b2001};
the requirement for the short time asymptotics of the constant coefficient equation in \cite{ep} is $0 <  \beta < 3 \sqrt{\alpha\gamma}$.}
\la{str conv}
\ee
In the recent article \cite{bb2018} sharp Gaussian estimates where obtained for the heat kernel of the operator (\ref{oper h}) without the strong convexity assumption.
Short time asymptotics were also obtained from which follows in particular that there is no absolute sharp exponential constant but instead the best constant depends on the range of the function
\be
Q(x)=\frac{\beta(x)}{\sqrt{\alpha(x)\gamma(x)}} \, , \qquad x\in\Omega \, .
\la{q(x)}
\ee

Our aim in the present article is to extend the estimates of \cite{bb2018} to the case where the operator $H$ is not uniformly elliptic and/or is not self-adjoint; in particular a sharp exponential constant
is obtained.
Concerning the singularity or degeneracy, we assume that $H$
is locally uniformly elliptic and that there is a positive weight function $w(x)$ that controls in a suitable sense the behaviour of the coefficients of the operator. Our main assumption consists of two general conditions (H1) and (H2) on $w(x)$, a weighted Sobolev inequality and a weighted interpolation inequality.
These conditions were introduced in \cite{b1998} in order to obtain (non-sharp) Gaussian estimates for non-uniformly elliptic self-adjoint operators. Besides conditions (H1) and (H2) we shall assume that the symbol $A(x,\xi)$ is close in an appropriate sense to a certain class of a ``good'' symbols induced by $w(x)$. These symbols correspond to operators which
additionally are self-adjoint and their coefficients are locally Lipschitz, with the behaviour near $\partial\Omega$ (or at infinity) being controlled by the weight $w(x)$. The estimates obtained herein
complement analogous estimates in \cite{b2004} where non-uniformly elliptic operators with strongly convex symbol were considered. The sharpness of the exponential constant $\sigma_*$ in our Gaussian estimate follows from the asymptotic estimates of  \cite{bb2018}.

The proof is based on Davies' exponential perturbation method. One has to consider three different regimes depending on the values taken
by the function $Q(x)$, namely $0\leq Q(x)\leq 3$ (the strongly convex regime), $Q(x)\leq 0$ and $Q(x)\geq 0$.
While the operator $H$ may be singular or degenerate, our assumptions guarantee that the function $Q(x)$ is bounded away from zero and infinity, which is crucial for the
implementation of the method.

\section{Heat kernel estimates}

\subsection{Setting and statement of main theorem}

Let $\Omega\subset\R^2$ be open and connected. We consider a differential operator $H$ on $L^2(\Omega)$ (complex-valued functions) given formally by
\be
Hu(x)=\partial_{x_1}^2\big(\alpha(x)\partial_{x_1}^2u\big)+2\partial^2_{x_1x_2}(\beta(x)\partial^2_{x_1x_2}u)+\partial_{x_2}^2\big(\gamma(x)\partial_{x_2}^2u\big),
\la{eq}
\ee
where $\alpha$, $\beta$ and $\gamma$ are complex-valued, locally bounded functions on $\Omega$. In case $\Omega\neq \R^2$ we impose Dirichlet boundary conditions on $\po$. 
The operator $H$ is defined by means of the quadratic form
\[
Q(u)=\int_{\Omega}\big\{\alpha(x)|u_{x_1x_1}|^2+2\beta(x)|u_{x_1x_2}|^2+\gamma(x)|u_{x_2x_2}|^2\big\}\,dx,
\]
defined initially on $C^{\infty}_c(\Omega)$. We assume that there exists a positive weight $w(x)$ with $w^{\pm 1}\in L^{\infty}_{loc}(\Omega)$ that controls the functions $\alpha(x),\beta(x),\gamma(x)$ in the following sense:
First, there holds
\be
|\alpha(x)|\leq cw(x),\qquad|\beta(x)|\leq cw(x),\qquad|\gamma(x)|\leq cw(x), \qquad x\in\Omega,
\la{5}
\ee
for some $c>0$ and second, the weighted G{\aa}rding inequality
\[
\re\,Q(u)\geq c\int_{\Omega}w(x)|\nabla^2u|^2\,dx,\qquad u\in C^{\infty}_c(\Omega),
\]
is valid for some $c>0$ (here $\nabla^2u$ denotes the vector whose components are the second-order partial derivatives of $u$).
This implies \cite[Theorem 7.12]{a} an analogous inequality for the symbol $A(x,\xi)$ of $H$, namely
\[
{\rm Re}\, A(x,\xi) \geq c \, w(x) |\xi|^{4} \; , \qquad x\in\Omega\, , \; \xi\in\R^2.
\]
The quadratic form $Q$ is closable and the domain of the closure is a weighted Sobolev space which we denote by $H^{2}_{w,0}(\Omega)$. We retain the same 
symbol, $Q$, for the closure of the above form and define
$H$ as the associated accretive operator on $L^2(\Omega)$, so that $\langle Hf,f\rangle=Q(u)$,
$f\in{\rm Dom}(H)$, and $Hu$ is given by (\ref{eq}) in the weak sense.

We make two assumptions on the weight $w(x)$, a weighted Sobolev inequality and a weighted interpolation inequality:

\

(H1)	There exist $s\in [\frac{1}{2}, 1]$ and $c>0$ such that
\[
\|u\|_{\infty}\leq c[\re Q(u)]^{\frac{s}{2}}\|u\|^{1-s}_{2},\qquad u\in C^{\infty}_c(\Omega).
\]

(H2)	There exists a constant $c>0$ such that
\[
\int_{\Omega}w^{\frac{1}{2}}|\nabla u|^2\,dx\, \leq
\,\epsilon\,\int_{\Omega}w|\nabla^{2}u|^2\,dx+c\epsilon^{-1}\,\int_{\Omega}|u|^{2}\,dx,
\]
for all $0<\epsilon<1$ and all $u\in C^{\infty}_c(\Omega)$.

Both (H1) and (H2) are satisfied when $H$ is uniformly elliptic, in which case the best value for the exponent $s$ is $s=1/2$, showing that in the general case we cannot expect any value that is smaller than $1/2$; in particular, (H1) is valid with $s=1/2$ if $w(x)$ is bounded away from zero. We refer to \cite{b1998} for a more detailed discussion of these conditions, including examples where they are both valid.

We note that condition (H2) implies that for any $k$, $l$ with $0\leq k$, $l\leq 2$, $k+l<4$, there exists a constant $c>0$ such that
\be
(1+\lambda^{4-k-l})\int_{\Omega}w^{\frac{k+l}{4}}|\nabla^{k}u|\,|\nabla^{l}u|\,dx\,
\leq \,\epsilon\,\re Q(u)+c\epsilon^{-\frac{k+l}{4-k-l}}(1+\lambda^{4})\|u\|^{2}_{2},
\la{eq:8}
\ee
for all $\epsilon\in (0,1)$, $\lambda>0$ and all $u\in C^{\infty}_c(\Omega)$. Indeed, for $\lambda=1$, (\ref{eq:8}) is a consequence of (H2) and the Cauchy-Schwarz inequality; the case $\lambda<1$ follows trivially from the case $\lambda=1$; finally, writing (\ref{eq:8}) for $\lambda=1$ and replacing $\epsilon$ by $\epsilon\lambda^{k+l-4}$ we obtain the result for $\lambda>1$.

We define the weighted Sobolev space
\[
W^{1,\infty}_{w}(\Omega)=\{u\in W^{1,\infty}_{{\rm loc}}(\Omega):    \exists c\geq 0  : \;  |u(x)|\leq c \, w(x) \; , \; |\nabla u(x)|\leq c\, w(x)^{\frac{3}{4}},\;x\in\Omega\}.
\]
{\bf Definition 1.} We say that the symbol $A(x,\xi)$ lies in $\mathcal{G}_w$ if the functions $\alpha(x)$, $\beta(x)$, $\gamma(x)$
are real-valued and belong in $W^{1,\infty}_{w}(\Omega)$.

We think of $\mathcal{G}_w$ as a class of ``good'' symbols. By assumption (\ref{5}) the last condition holds true if and only if
\[
|\nabla\alpha(x)|+|\nabla\beta(x)|+|\nabla\gamma(x)|\leq cw(x)^{\frac{3}{4}} \; , \quad x\in\Omega \, .
\]
To state our main  result we need some more definitions. We first set
\[
\cE_{w} =\big\{\phi\in C^{2}(\Omega)\cap L^{\infty}(\Omega):\; \phi \mbox{ real valued, }  \exists c>0   : \,
|\nabla\phi| \leq  c w^{-\frac{1}{4}}  ,  \;   |\nabla^2\phi| \leq  c w^{-\frac{1}{2}} \big\}.
\]
In case where the symbol $A(x,\xi)$ belongs in $\mathcal{G}_w$ (so in particular it is real-valued) we additionally define for any $M>0$  the subclass
\[
\cE_{A,M} =\big\{\phi\in \cE_{w} \, :\;  A(x,\nabla\phi(x))\leq 1,\,|\nabla^{2}\phi(x)|\leq M\,w(x)^{-\frac{1}{2}}, \; x\in\Omega\big\} \, ;
\]
our Gaussian estimates will be expressed in terms of the distance
\[
d_M(x,x')=\sup\big\{\phi(x')-\phi(x)\,:\;\;\phi\in\cE_{A,M}\big\}
\]
for arbitrariy large (but finite) $M$; we note that as $M\to +\infty$ this converges to
\[
d(x,x') = \sup \{ \phi(x') -\phi(x) \; : \;   \phi\in {\rm Lip}(\Omega) \, , \;\;  A( y ,\nabla\phi(y))\leq 1 \, , \;\; {\rm a.e. }\;\; y\in\Omega\}.
\]
The domain $\Omega$ is essentially partitioned in three components depending on the values of the bounded function
$Q(x)$ (cf. (\ref{q(x)})).
In particular, assuming always that the symbol $A(x,\xi)$ belongs in $\mathcal{G}_w$, we define the locally Lipschitz functions
\[
k(x)=\tarr{8\frac{1-Q(x)}{(1+Q(x))^2},}{\mbox{ if } Q(x) \leq 0,}{8,}{\mbox{ if }0\leq Q(x)\leq3,}{Q(x)^2-1,}{\mbox{ if }Q(x)\geq 3,}
\]
and
\[
\sigma(x)=\frac{3}{4}\cdot\Big(\frac{1}{4k(x)}\Big)^{1/3}= \tarr{ \frac{3}{8\cdot 4^{1/3}} \frac{ (1+Q(x))^{2/3}}{ (1-Q(x))^{1/3}},} 
{\mbox{ if }Q(x) \leq 0,}{\frac{3}{8\cdot 4^{1/3}},}{\mbox{ if }0\leq Q(x)\leq3,}{ \frac{3}{4^{4/3}} (Q(x)^2-1)^{-1/3} ,}{\mbox{ if }Q(x) \geq 3.}
\]
We also set 
\[
k^*=\sup_{x\in\Omega}k(x) \quad \mbox{ and }\quad
\sigma_*= \inf_{x\in\Omega}\sigma(x) = \frac{3}{4}\cdot\Big(\frac{1}{4k^*}\Big)^{1/3}.
\]
In the general case where the symbol does not belong in $\mathcal{G}_w$ we denote by $\theta$ the
following weighted distance of the symbol $A(x,\xi)$ from $\mathcal{G}_w$,
\[
\theta=\inf_{\tilde{A} \in \mathcal{G}_w} \sup_{\Omega} \, \max_{|\xi|=1}   \frac{ \big|   A(x,\xi) -  \tilde{A}(x,\xi)   \big|}{w(x)} \, .
\]
We shall think of $\theta$ as a small number.

We now state our main result; the constants $c_{\epsilon}$, $c_{\epsilon,M}$ may also depend on the operator $H$.
\begin{theorem}
\label{thm1}
Assume that (H1) and (H2) are satisfied. \nl
(a) Assume that the symbol $A(x,\xi)$ belongs in $\mathcal{G}_w$. Then for all $\epsilon\in(0,1)$ and all $M$ large there exist $c_\epsilon,c_{\epsilon,M}<\infty$ such that
\be
|G(x,x',t)|\leq c_\epsilon t^{-s}\exp\Big\{-(\sigma_*-\epsilon)\frac{d_M(x,x')^{4/3}}{t^{1/3}}+c_{\epsilon,M}t\Big\},
\la{cov1}
\ee
for all $x,x'\in\Omega$ and $t>0$. \nl
(b) If $A(x,\xi)$ does not belong $\mathcal{G}_w$ then there exists $c>0$ such that for all $\epsilon\in(0,1)$ 
and all $M$ large there exist $c_\epsilon,c_{\epsilon,M}<\infty$ such that
\[
|G(x,x',t)|\leq c_\epsilon t^{-s}\exp\Big\{-(\sigma_*-c\theta-\epsilon)
\frac{d_M(x,x')^{4/3}}{t^{1/3}}+c_{\epsilon,M}t\Big\},
\]
for all $x,x'\in\Omega$ and $t>0$; here $\sigma_*$ and $d_M(x,x')$ are defined as above corresponding to a symbol $\tilde{A}(x,\xi)$ in $\mathcal{G}_w$
for which $|A(x,\xi) - \tilde{A}(x,\xi)| \leq 2\theta w(x)|\xi|^{4}$, $x\in\Omega$, $\xi\in\R^2$.
\end{theorem}

{\em Remarks.} (1) It follows from the asymptotic estimates obtained in \cite{bb2018} that the constant $\sigma_*$ is the best possible. \nl
(2) In case (b) one could define the exponential constant $\sigma_*$  and the distance $d_M(x,x')$ using the symbol $A(x,\xi)$ rather than 
$\tilde{A}(x,\xi)$. The resulting estimate would be comparable to the one in the theorem; such differences are anyway absorbed in the term $c\theta$
in the exponential and we prefer to used $\tilde{A}(x,\xi)$ for the definition of these quantities since otherwise the proofs would be longer.


\subsection{Proof of Theorem \ref{thm1}}

As already mentioned, the proof makes use of Davies' perturbative argument \cite{d}.
It follows from hypothesis (H2) that for any  $\psi\in\cE_{w}$ the (multiplication) operator $e^{\psi}$ leaves the Sobolev space
$H^{2}_{w,0}(\Omega)$ invariant so we may define a sesquilinear form $Q_{\psi}$ on $H^{2}_{w,0}(\Omega)$ by
$Q_{\psi}(u)=Q(e^{\psi}u,e^{-\psi}u) $;
here
\[
Q(u,v)= \int_{\Omega}\big\{\alpha(x)u_{x_1x_1}\overline{v}_{x_1x_1}+2\beta(x)u_{x_1x_2}\overline{v}_{x_1x_2}
+\gamma(x)u_{x_2x_2}\overline{v}_{x_2x_2}\big\}\,dx
\]
is the sesquilinear form associated to $Q(\cdot)$, hence
\bea
Q_{\psi}(u)&=&\int_{\Omega}\Big[\alpha(x)(e^{\psi}u)_{x_1x_1}
(e^{-\psi}\overline{u})_{x_1x_1}+2\beta(x)(e^{\psi}u)_{x_1x_2}(e^{-\psi}\overline{u})_{x_1x_2}\nonumber\\
&&\quad\quad+\gamma(x)(e^{\psi}u)_{x_2x_2}(e^{-\psi}\overline{u})_{x_2x_2} \Big] \,dx.  \la{eenndd}
\eea
We shall need the following result, see \cite[Proposition 3.2]{b2004}:
\begin{lemma}
Assume that (H1) and (H2) are satisfied. Let $\psi\in\cE_w$ be given and let $k \in\R$ be such that
\[
{\rm Re}\,Q_{\psi}(u)\geq - k  \,\|u\|_2^2
\]
for all $u\in C_c^{\infty}(\Omega)$. Then for any $\delta\in(0,1)$ there exists a constant $c_\delta$ such that
\[
|G(x,x',t)|\leq c_{\delta} t^{-s}\exp\big\{\psi(x)-\psi(x')+(1+\delta) k t\big\},
\]
for all $x,x'\in\Omega$ and all $t>0$.
\label{lem:ebd}
\end{lemma}
We now take in (\ref{eenndd})  $\psi =\lambda\phi$ where $\lambda>0$ and $\phi\in\cE_{A,M}$.
After expanding, the exponentials $e^{\lambda\phi}$ and $e^{-\lambda\phi}$ cancel and we obtain that $Q_{\lambda\phi}(u)$ is a linear combination of terms of the form
\be
\lambda^s \int_{\Omega}b_{s,\gamma,\delta}(x)D^{\gamma}u\,D^{\delta}\overline{u}\,dx,
\la{eq:19}
\ee
(multi-index notation) where $s+|\gamma+\delta|\leq 4$ and
each function $b_{s\gamma\delta}(x)$ is a product of one of the functions $\alpha(x)$, $\beta(x)$, $\gamma(x)$ and first or second order derivatives of $\phi(x)$ (see
also (\ref{789}) below). Recalling (\ref{5}) we see that for each such term we have
\be
|b_{s,\gamma,\delta}(x) | \leq c w(x)^{\frac{|\gamma+\delta|}{4}} \; , \qquad x\in\Omega \, .
\la{est}
\ee

\noindent
{\bf Definition 2.} We denote by $\cL$ the space of (finite) linear combinations of terms of the form (\ref{eq:19})
with $s+|\gamma+\delta|<4$ and $|b_{s,\gamma,\delta}(x)|\leq cw(x)^{\frac{|\gamma+\delta|}{4}}$.

We note that if the form (\ref{eq:19}) belongs in $\cL$, then by (\ref{eq:8})  we have for any $\epsilon>0$,
\begin{align}
|T(u)| & \leq  c\lambda^s \int_{\Omega}  cw(x)^{\frac{|\gamma+\delta|}{4}}    |D^{\gamma}u|  \, |D^{\delta}u| dx \nonumber \\
&\leq    \epsilon\,\re Q(u)+c\epsilon^{-\frac{|\gamma +\delta|}{4-|\gamma +\delta|}}(1+\lambda^{\frac{4s}{4-|\gamma +\delta|}})  \|u\|^{2}_{2} \nonumber \\
&\leq    \epsilon\,\re Q(u)+c\epsilon^{-3}(1+\lambda^3)  \|u\|^{2}_{2}  \; .
\la{form:l}
\end{align}

We next define the quadratic form
\begin{align}
\la{789}&\qquad Q_{1,\lf}(u)= \nonumber \\
&\int_{\Omega}\Big\{\lambda^4\big[\alpha(x)\phi_{x_1}^4+2\beta(x)\phi_{x_1}^2\phi_{x_2}^2+\gamma(x)\phi_{x_2}^4\big]|u|^2\nonumber \\
&\quad+\lambda^2\Big\{\alpha(x)\phi_{x_1}^2(u\overline{u}_{x_1x_1}+u_{x_1x_1}\overline{u} -4|u_{x_1}|^2)\nonumber \\
&\quad+2\beta(x)\big[\phi_{x_1}\phi_{x_2}(u\overline{u}_{x_1x_2}+u_{x_1x_2}\overline{u}-u_{x_1}\overline{u}_{x_2}-u_{x_2}\overline{u}_{x_1})-(\phi_{x_2}^2|u_{x_1}|^2
+\phi_{x_1}^2|u_{x_2}|^2)\big]\nonumber\\
&\quad+\gamma(x)\phi_{x_2}^2(u\overline{u}_{x_2x_2}+u_{x_2x_2}\overline{u}-4|u_{x_2}|^2)\Big\}\nonumber \\
&\quad +\alpha(x)|u_{x_1x_1}|^2+2\beta(x)|u_{x_1x_2}|^2+\gamma(x)|u_{x_2x_2}|^2\Big\}\,dx. 
\end{align}
It may be seen that $Q_{1,\lf}(u)$ contains precisely those terms of the form (\ref{eq:19}) from the expansion of $Q_{\lf}(u)$ for which we have $s+|\gamma+\delta|=4$.
Hence, recalling also (\ref{est}), the difference $Q_{\lf}(\cdot)-Q_{1,\lf}(\cdot)$ belongs  in $\cL$.

We now define the polar symbol
\[
A(x,z,z')=\alpha(x)z_1^2z_1'^2+2\beta(x)z_1z_2z_1'z_2'+\gamma(x)z_2^2z_2'^2,\qquad x\in\Omega,\quad z,\,z'\in\C^2 \, .
\]
We note that for $z=z'=\xi\in\R^2$ this reduces to the symbol $A(x,\xi)$ of $H$.
For $x\in\Omega$ and $\xi,\xi',\eta\in\R^2$ we also set
\be
S(x,\xi,\xi',\eta)={\rm Re}\,A(x,\xi+i\eta,\xi'+i\eta)+k(x)A(x,\eta).
\la{star}
\ee
Given $\phi\in\cE_w$ and $\lambda >0$ we define the quadratic form $S_{\lf}$ on $H^{2}_{w,0}(\Omega)$ by
\[
S_{\lf}(u)=\frac{1}{(2\pi)^{2}}\iiint_{\Omega\times\R^2\times\R^2}S(x,\xi,\xi',\lambda\nabla\phi)e^{i(\xi-\xi')\cdot x}\hat{u}(\xi)\overline{\hat{u}(\xi')}\,d\xi\,d\xi'\,dx.
\]
\begin{lemma}
\label{lem2}
Assume that the symbol $A(x,\xi)$ lies in $\mathcal{G}_w$.
Let $\phi\in\cE_w$ and $\lambda>0$. There holds
\[
{\rm Re}\, Q_{1,\lf}(u)+\int_{\Omega} k(x)A(x,\lambda\nabla\phi)|u|^2dx=S_{\lf}(u),
\]
for all $u\in C^{\infty}_c(\Omega)$.
\end{lemma}
\noindent
{\em Proof.} This follows from (\ref{star}) by using the relation
\[
D^{\alpha}u(x) =(2\pi)^{-1}\int_{\R^2}(i\xi)^{\alpha}e^{ix\cdot\xi}\hat{u}(\xi)d\xi
\]
for the various terms that appear in $Q_{1,\lf}$; the fact that $\alpha(x)$, $\beta(x)$ and $\gamma(x)$ are real-valued is also used here. $\hfill\Box$

We now define for each $x\in\Omega$ a quadratic form $\Gamma(x , \cdot)$ in $\C^6$ by
\bean
&&\Gamma(x,p)=\\
&&\hspace{-0.5cm}\left\{\begin{array}{l}{(Q+1)|p_1|^2+(Q+1)|p_2|^2-Q|p_3|^2-2Q|p_4|^2-}\\[0.1cm]
{\hspace{3cm}-2Q|p_5|^2-\frac{Q(3-Q)^2}{(1+Q)^2}|p_6|^2,
\hspace{1.6cm}\mbox{if }-1<Q(x)<0,}\\[0.3cm]
\frac{3-Q}{3}|p_1|^2+\frac{3-Q}{3}|p_2|^2+\frac{Q}{3}|p_1+p_2|^2+\frac{4Q}{3}|p_3|^2,
\hspace{2cm}\mbox{if }0\leq Q(x)\leq3, \\[0.3cm]
{2(Q-3)|p_1|^2+|p_2|^2 +2(Q-1)|p_3|^2+2\frac{Q-3}{Q-1}(Q+1)(Q^2+3)|p_4|^2,}\\[0.1cm]
\hspace{9cm}\mbox{if }Q(x)>3,
\end{array}
\right.
\eean
for any $p=(p_1,\ldots,p_6)\in\C^6$.
Clearly $\Gamma(x , \cdot)$ is positive semidefinite for each $x\in\Omega$. We denote by $\Gamma(x,\cdot,\cdot)$ the corresponding sesquilinear form in $\C^6$, that is $\Gamma(x, p ,q)$ is given by a formula similar to the one above with each $|p_k|^2$ being replaced by $p_k\overline{q_k}$
and with $|p_1+p_2|^2$ being replaced by $(p_1+p_2)\overline{(q_1+q_2)}$.

Next, for any $x\in\Omega$ and $\xi,\eta\in\R^2$ we define a vector $p_{x,\xi,\eta}\in\R^6$ by
\bean
&&p_{x,\xi,\eta}=\\
&& \left\{
\begin{array}{l}
{\!\!\!\!\Big(\alpha^{1/2}[\xi_1^2-\frac{3-Q}{1+Q}\eta_1^2],\gamma^{1/2}[\xi_2^2-\frac{3-Q}{1+Q}\eta_2^2],\,\alpha^{1/2}\xi_1^2-\gamma^{1/2}\xi_2^2,\,\alpha^{1/2}\xi_1\eta_1 \! +\!\gamma^{1/2}\xi_2\eta_2,}\\
{\qquad \alpha^{1/4}\gamma^{1/4}(\xi_1\eta_2+\xi_2\eta_1),\,\alpha^{1/2}\eta_1^2-\gamma^{1/2}\eta_2^2\Big),}
\hspace{1.3cm}{\mbox{ if } -1<Q(x)<0,} \\[0.2cm]
{\!\!\!\Big(\alpha^{1/2}[\xi_1^2-3\eta_1^2],\,\gamma^{1/2}[\xi_2^2-3\eta_2^2],\,\alpha^{1/4}\gamma^{1/4}[\xi_1\xi_2-3\eta_1\eta_2],\,0,\,0,\,0\Big),} \\
\hspace{8.3cm}{\mbox{ if } 0\leq Q(x)\leq3,} \\[0.2cm] 
{\!\!\!\Big(\alpha^{1/2}\xi_1\eta_1-\gamma^{1/2}\xi_2\eta_2,\,\alpha^{1/2}(\xi_1^2-Q \eta_1^2)+\gamma^{1/2}(\xi_2^2-Q\eta_2^2),}\\
{\qquad\alpha^{1/4}\gamma^{1/4}[\xi_1\xi_2-\frac{Q+3}{Q-1}\eta_1\eta_2],\,\alpha^{1/4}\gamma^{1/4}\eta_1\eta_2,\,0,\,0\Big),}
\hspace{.45cm} {\mbox{ if }Q(x)>3.}
\end{array}
\right.
\eean
A crucial property of the form $\Gamma(x,\cdot)$ and the vectors $p_{x,\xi,\eta}$ is that
\be
S(x;\xi,\xi,\eta)=\Gamma(x,p_{x,\xi,\eta},p_{x,\xi,\eta}),
\label{s:g}
\ee
for all $x\in\Omega$ and $\xi,\eta\in\R^2$.

We finally define a quadratic form $\Gamma_{\lf}(\cdot)$ on $H^{2}_{w,0}(\Omega)$ by
\[
\Gamma_{\lf}(u)=\frac{1}{(2\pi)^{2}}\iiint_{\Omega\times\R^2\times\R^2}\Gamma(x, \, p_{x,\xi,\lambda\nabla\phi},p_{x,\xi',\lambda\nabla\phi})e^{i(\xi-\xi')\cdot x}\hat{u}(\xi)
\overline{\hat{u}(\xi')}\,d\xi\,d\xi'\,dx.
\]
We then have
\begin{lemma}
\label{lem3}
Assume that the symbol $A(x,\xi)$ lies in $\mathcal{G}_w$. Then the difference $S_{\lf}(\cdot)-\Gamma_{\lf}(\cdot)$ belongs to $\cL$.
\end{lemma}
\noindent
{\em Proof.} We consider the difference
\[
S(x,\xi,\xi',\eta)-\Gamma(x,p_{x,\xi,\eta},p_{x,\xi',\eta}),
\]
of the two symbols and we group together terms that have the property that if we set $\xi'=\xi$ then they are similar as monomials of the variables
$\xi$ and $\eta$. Due to (\ref{s:g}) one can use integration by parts to conclude that the total contribution of each such
group belongs to $\cL$. We shall illustrate this for two particular groups, the one consisting of terms which for $\xi=\xi'$ involve the monomial $\xi_1^2\eta_1^2$ and those which
for $\xi=\xi'$ involve $\xi_1^2\eta_2^2$. For the sake of brevity we shall consider directly the sum of the terms of both groups.

The terms of these two groups from $S(x,\xi,\xi',\eta)$ add up to
\[
-\alpha(x)\eta_1^2(\xi_1^2+\xi_1'^2+4\xi_1\xi_1')-2\beta(x)\eta_2^2\xi_1\xi_1'.
\] 

The corresponding terms in $\Gamma(x,p_{x,\xi,\eta},p_{x,\xi',\eta})$ are
\[
\left\{
\begin{array}{l}
{\!\!\alpha(x)\eta_1^2\big[(Q(x)-3)(\xi_1^2+\xi_1'^2)-2Q(x)\xi_1\xi_1'\big]-2\beta(x)\eta_2^2\xi_1\xi_1',} \hspace{1.4cm} {\mbox{if }Q(x)\leq 0,} \\[0.2cm]
{\!\!-3\alpha(x)\eta_1^2(\xi_1^2+\xi_1'^2)-\beta(x)\eta_2^2(\xi_1^2+\xi_1'^2),}
\hspace{3.6cm}\mbox{if }0\leq Q(x)\leq3,\\[0.2cm]
{\!\!\alpha(x)\eta_1^2\big[-Q(x)(\xi_1^2+\xi_1'^2)+2(Q(x)-3)\xi_1\xi_1'\big]-\beta(x)\eta_2^2(\xi_1^2+\xi_1'^2),} \hspace{.3cm}\mbox{if }Q(x)\geq 3.
\end{array}
\right.
\]
Hence the difference of these terms in $S(x,\xi,\xi',\eta)-\Gamma(x,p_{x,\xi,\eta},p_{x,\xi',\eta})$ is
\[
\tarr{\alpha(x)\eta_1^2\big[2-Q(x)\big](\xi_1-\xi_1')^2,}{\mbox{ if }Q(x)\leq 0,}
{\big[2\alpha(x)\eta_1^2+\beta(x)\eta_2^2\big](\xi_1-\xi_1')^2,}{\mbox{ if }0\leq Q(x)\leq3,}
{\big[\alpha(x)(Q(x)-1)\eta_1^2+\beta(x)\eta_2^2\big](\xi_1-\xi_1')^2,}{\mbox{ if } Q(x)\geq 3.}
\]

This can also be written as $\big[\alpha(x)\eta_1^2 R(x)+\eta_2^2 P(x)\big](\xi_1-\xi_1')^2$ where
\[
R(x)=\tarr{2-Q(x),}{\mbox{ if }Q(x)\leq 0,}{2,}{\mbox{ if }0\leq Q(x)\leq3,}{Q(x)-1,}{\mbox{ if }Q(x)\geq 3,}
\quad \mbox{and} \quad
P(x)=\left\{\begin{array}{rcl} 0, & \mbox{if} & \beta(x) \leq 0, \\ \beta(x), &  \mbox{if} & \beta(x)\geq 0 . \end{array}\right.
\]

Inserting this in the triple integral and recalling that $\eta=\lambda\nabla\phi$ we obtain that the contribution of the above terms in the difference $S_{\lf}(u)-\Gamma_{\lf}(u)$ is
\bean
&&(2\pi)^{-2}\iiint_{\Omega\times\R^2\times\R^2}\big[\alpha(x)R(x)\phi_{x_1}^2+P(x)\phi_{x_2}^2\big](\xi_1-\xi_1')^2\lambda^2 e^{i(\xi-\xi')\cdot x} \\
&&\hspace{4cm}\hat{u}(\xi)\overline{\hat{u}(\xi')}\,d\xi\,d\xi'\,dx\\
&=&\lambda^2\int_{\Omega}\big[\alpha(x)R(x)\phi_{x_1}^2+P(x)\phi_{x_2}^2\big](-u_{x_1x_1}\overline{u}-u\overline{u}_{x_1x_1}-2|u_{x_1}|^2)dx\\
&=&-\lambda^2\int_{\Omega}\big[\alpha(x)R(x)\phi_{x_1}^2+P(x)\phi_{x_2}^2\big]( u_{x_1}\overline{u}+u\overline{u}_{x_1})_{x_1}dx \\
&=&\lambda^2\int_{\Omega}\big[\alpha(x)R(x)\phi_{x_1}^2+P(x)\phi_{x_2}^2\big]_{x_1}( u_{x_1}\overline{u}+u\overline{u}_{x_1})dx  ,
\eean
where we have used the fact that
the function $\alpha(x)R(x)\phi_{x_1}^2+P(x)\phi_{x_2}^2$ is locally Lipschitz. To conclude that the last expression belongs in  $\cL$ we must prove that (\ref{est}) is valid, that is
$\big|   [\alpha(x)R(x)\phi_{x_1}^2+P(x)\phi_{x_2}^2]_{x_1}   \big| \leq c w(x)^{\frac{1}{4}}$.
We shall only consider the first of the two terms, the proof being similar for the second. Using the relations
$|Q(x)|\leq c$, $|\nabla Q(x)|\leq cw(x)^{-1/4}$ we obtain
\begin{eqnarray*}
\big|  (\alpha(x)R(x)\phi_{x_1}^2)_{x_1} \big| &\leq &  |\alpha_{x_1} R|  \phi_{x_1}^2 + |\alpha R_{x_1} | \phi_{x_1}^2 + 2| \alpha R  \phi_{x_1} \phi_{x_1x_1}  | \\
&\leq & c w^{\frac{3}{4}} w^{-\frac{1}{2}} + cw w^{-\frac{1}{4}}  w^{-\frac{1}{2}}  + cw w^{-\frac{1}{4}} M w^{-\frac{1}{2}} \\
&=& c_M w^{\frac{1}{4}},
\end{eqnarray*}
as required.  $\hfill\Box$
\begin{lemma}
\label{lem4}
Assume that the symbol $A(x,\xi)$ lies in $\mathcal{G}_w$ and let $M>0$ be given. Then for any $\phi\in\cE_{A,M}$ and $\lambda>0$ we have
\[
\re\,Q_{\lambda\phi}(u)\geq-k^*\lambda^4\,\|u\|_2^2+T(u),
\]
for some quadratic form $T\in\cL$ and all $u\in C^{\infty}_c(\Omega)$.
\end{lemma}
\noindent
{\em Proof.} The assumption $\phi\in\cE_{A,M}$ implies that $A(x,\nabla\phi(x))\leq1$, $x\in\Omega$. Recalling that the difference $Q_{\lf}(\cdot)-Q_{1,\lf}(\cdot)$ belongs in $\cL$ and
using Lemmas \ref{lem2},  \ref{lem3} and \ref{lem4}
we obtain
\bean
\re\,Q_{\lambda\phi}(u)&=&-\int_\Omega k(x)A(x,\lambda\nabla\phi)\,|u|^2\,dx+\Gamma_{\lf}(u)+T(u)\\
&\geq&-k^*\lambda^4\int_\Omega|u|^2\,dx+\Gamma_{\lf}(u)+T(u),
\eean
for some form $T\in\cL$ and all $u\in C^{\infty}_c(\Omega)$. Moreover
\begin{align*}
\Gamma_{\lf}(u)\, =\; &\frac{1}{(2\pi)^{2}}\iiint_{\Omega\times\R^2\times\R^2}\Gamma(x, \, p_{x,\xi,\lambda\nabla\phi},p_{x,\xi',\lambda\nabla\phi})e^{i(\xi-\xi')\cdot x}\hat{u}(\xi)\overline{\hat{u}(\xi')}\,d\xi\,d\xi'\,dx\\
\, =\;&\frac{1}{(2\pi)^{2}}\int_{\Omega}\Gamma\Big(x , \,\int_{\R^2}e^{i\xi\cdot x}\hat{u}(\xi)p_{x,\xi,\lambda\nabla\phi}d\xi,\int_{\R^2}e^{i\xi'\cdot x}\hat{u}(\xi')p_{x,\xi',\lambda\nabla\phi}d\xi'\Big)\,dx\\
\, \geq \; &0,
\end{align*}
by the positive semi-definiteness of $\Gamma$; the result follows.$\hfill\Box$

\medskip

\noindent
{\bf\em Proof of Theorem \ref{thm1}.} 
{\em Part (a).} We claim that for any $\epsilon$ and $M$ positive there exists $c_{\epsilon,M}$ (which may also depend on the operator $H$) such that
\be
\re\,Q_{\lf}(u)\geq-\Big\{(k^*+\epsilon)\lambda^4+c_{\epsilon,M}(1+\lambda^3)\Big\}\|u\|_2^2.
\la{claim}
\ee
for all $\lambda>0$ and $\phi\in\cE_{A,M}$.
To prove this we first recall (cf.  (\ref{form:l})) that any form $T\in\cL$ satisfies
\[
|T(u)|\leq\epsilon Q(u)+c_{\epsilon,M}(1+\lambda^3)\,\|u\|_2^2,
\]
for all $\epsilon\in(0,1)$, $\lambda>0$ and $u\in C^{\infty}_c(\Omega)$. Hence, since $Q(u)$ is real, Lemma \ref{lem4}  implies
\be
\la{gui1}
\re\,Q_{\lf}(u)\geq-\Big\{k^*\lambda^4+c_{\epsilon,M}(1+\lambda^3)\Big\}\|u\|_2^2-
\epsilon Q(u).
\ee
Now, considering the expansion of $Q_{\lf}$ already discussed and recalling (\ref{eq:8}) we infer that there exists a constant $c_M$ such that for any $\phi\in\cE_{A,M}$ and $\lambda>0$ there holds
\be
\big| Q(u)-Q_{\lf}(u) \big| \leq\frac{1}{2}Q(u)+c_M(\lambda+\lambda^4)\|u\|_2^2 \, .
\la{fil1}
\ee
Furthermore, we note that the dependence on $M$ in this estimate comes from those terms in the expansion 
of $Q_{\lf}$ that contain at least one second-order derivative of $\phi$.
Since the coefficient of $\lambda^4$ in the expansion only involves first 
derivatives  of $\phi$, (\ref{fil1}) can be improved to
\[
\big|Q(u)-Q_{\lf}(u) \big|\leq\frac{1}{2}Q(u)+\big\{c_M(\lambda+\lambda^3)+c
\lambda^4\big\}\|u\|_2^2,
\]
which in turn implies
\be
Q(u)\leq\,  2{\rm Re}\,Q_{\lf}(u)+\big\{c_M(\lambda+\lambda^3)+c\lambda^4\big\}
\|u\|_2^2.
\la{last}
\ee
Let $u\in C^{\infty}_c(\Omega)$ be given. If $\re Q_{\lf}(u)\geq0$ then (\ref{claim}) is obviously true. If not we then have from (\ref{gui1}) and (\ref{last})
\bean
\re\,Q_{\lf}(u)&\geq&-\Big\{k^*\lambda^4+c_{\epsilon,M}(1+\lambda^3)\Big\}\|u\|_2^2-2\epsilon\, \re Q_{\lf}(u) \\
&& -\epsilon\big\{c_M(\lambda+\lambda^3)+c\lambda^4\big\}\|u\|_2^2\\
&\geq&-\Big\{(k^*+c\epsilon)\lambda^4+c_{\epsilon,M}(1+\lambda^3)+\epsilon\big\{c_M(\lambda+\lambda^3)+c\lambda^4\big\}\Big\}\|u\|_2^2,
\eean
and (\ref{claim}) again follows; hence the claim has been proved.

We complete the standard argument; Lemma \ref{lem:ebd} and (\ref{claim}) imply
\[
|G(x,x',t)|<c_\epsilon t^{-s}\exp\Big\{\lambda \big( \phi(x)-\phi(x') \big)+(1+\epsilon)\big\{(k^*+\epsilon)\lambda^4+c_{\epsilon,M}(1+\lambda^3)\big\}\,t\Big\}, 
\]
for all $\epsilon\in(0,1)$. Optimizing over $\phi\in\cE_{A,M}$ yields
\[
|G(x,x',t)|<c_\epsilon t^{-s}\exp\Big\{-\lambda d_M(x,x')+(1+\epsilon) \big\{ (k^*+\epsilon)\lambda^4+c_{\epsilon,M}(1+\lambda^3) \big\}\,t\Big\} .
\]
Finally choosing $\lambda=[d_M(x,x')/(4k^*t)]^{1/3}$ we have
\[
-\lambda d_M(x,x')+k^*\lambda^4t =-\sigma_*\frac{d_M(x,x')^{4/3}}{t^{1/3}},
\]
and (\ref{cov1}) follows.

{\em Part (b).} There exists a symbol $\tilde{A}(x,\xi)$ in  $\mathcal{G}_w$ such that
\[
 \max \big\{ | \alpha(x)-\tilde\alpha(x)| \, , \; |\beta(x)-\tilde\beta(x)|\,  , \; |\gamma(x)-\tilde\gamma(x)| \big\}  \leq 2\theta \, w(x) \; , \qquad x\in\Omega.
\]
Given $\phi\in\cE_{\tilde{A},M}$ and $\lambda>0$ it follows from the proof of Part (a) that
\be
\re\,\tilde{Q}_{\lf}(u)\geq-\Big\{ k^* \lambda^4+c_{\epsilon,M}(1+\lambda^3)\Big\}\|u\|_2^2-\epsilon \, \re Q(u),
\la{25}
\ee
for all $u\in C^{\infty}_c(\Omega)$. Moreover it is easily seen that
\be
 \big| Q_{\lf}(u)-\tilde{Q}_{\lf}(u) \big| \leq c\theta \big\{ \re Q(u)+\lambda^4\|u\|_2^2\big\}.
\la{26}
\ee
The argument used for (\ref{last}) also applies to $H$ and we thus obtain
\be
{\rm Re} \, Q(u)\leq\,  2{\rm Re}\,Q_{\lf}(u)+\big\{c_M(\lambda+\lambda^3)+c\lambda^4\big\} \|u\|_2^2.
\la{last1}
\ee
Combining (\ref{25}), (\ref{26}) and (\ref{last1}) we conclude that
\[
\re\,Q_{\lf}(u)\geq-\Big\{(k^*+c\theta+\epsilon )\lambda^4+c_{\epsilon,M}(1+\lambda^3)\Big\}\|u\|_2^2 , \qquad u\in C^{\infty}_c(\Omega),
\]
and the argument is completed as in Part (a); we omit further details. $\hfill\Box$

\end{document}